\documentclass{GANITA}

\citesort
\usepackage{subcaption}

\theoremstyle{cupthm}
\newtheorem{thm}{Theorem}[section]

\newtheorem{cor}[thm]{Corollary}

\theoremstyle{cupdefn}
\newtheorem{defn}[thm]{Definition}
\theoremstyle{cuprem}
\newtheorem{rem}[thm]{Remark}
\numberwithin{equation}{section}

\begin{document}
\runningtitle{A Bivariate Polynomial Problem for Matrices}
\title{A bivariate polynomial problem for matrices}
\author[1]{Dharm Prakash Singh}
\address[1]{Department of Applied Mathematics, School of Vocational Studies and Applied Sciences, Gautam Buddha University, Greater Noida, Uttar Pradesh,India-201312\email{prakashdharm@gmail.com}}
\author[2]{Amit Ujlayan$^*$}
\address[2]{Department of Applied Mathematics, School of Vocational Studies and Applied Sciences, Gautam Buddha University, Greater Noida, Uttar Pradesh,India-201312\email{iitk.amit@gmail.com}}
\author[3]{Bhim Sen Choudhary}
\address[3]{Department of Mathematics \& Astronomy, University of Lucknow, Lucknow, India-226007\email{bhim2051@gmail.com}}

\authorheadline{D. P. Singh, A. Ujlayan and B. S. Choudhary}


\support{$^*$ Corresponding Author\\Dharm Prakash Singh and Bhim Sen Choudhary thank University Grants Commission (UGC), India, for paying for research scholarship under NET/JRF scheme with Ref. Nos.: 1086/(OBC)(CSIR-UGC NET DEC. 2016) and 1251/(OBC)(CSIR-UGC NET DEC. 2016), respectively}

\begin{abstract}
This article proposes a bivariate polynomial problem for finite-order real matrices that endows a \textit{`sufficient condition'} for a map from the standard vector spaces of finite-order real matrices to the same dimensional bivariate polynomial subspaces (BVPSs) to be an isomorphism in some finite-dimensional BVPSs. In the process of solving, the article deals with the existence, uniqueness, and construction of the polynomials in some finite-dimensional BVPSs concerning the solution of the proposed problem. To this end, a relationship is established between the proposed problem and a class of Lagrange bivariate polynomial interpolation problems (LBVPIPs). As a result, the existence of a standard and a new class of finite-dimensional BVPSs of various total degrees has been established in which the proposed problem always possesses a unique solution. In addition, some formulas are derived to construct the needed polynomials in these BVPSs. Further, the possible applicability of the proposed problem is discussed in LBVPIPs, focusing on the finite rectangular schemes of bivariate interpolation points on the natural Cartesian grid. At last, some numerical examples are considered to justify the theoretical findings. 
\end{abstract}


\maketitle

\section{Introduction}\label{intro}
Over the field of real numbers $\mathcal{R}$, let $\Pi^{2}$ be the space of all bivariate polynomials, $\Pi_{k}^{2}$ be $(k+1)(k+2)/2$-dimensional subspace of $\Pi^{2}$ of degree at most $k$, $\Pi_{k}$ be $(k+1)$-dimensional vector space of all univariate polynomials of degree at most $k$, and $\mathcal{R}^{m \times n}$ be $mn$-dimensional vector space of all $m$ by $n$ matrices, where $m,n,k \in \mathbb{N}$ (set of all positive integers), \cite{1,2}. 

The bivariate polynomials are used to improve mathematical modeling in medical applications \cite{3}, image processing \cite{4,5}, signal processing \cite{5}, computer graphics \cite{6}, computer-aided design in control systems \cite{7,8}, etc., in which ``data" is stored as 2D arrays or finite-order real matrices. Therefore, by constructing corresponding bivariate polynomials that can map the data points to mathematical space, one can conduct advanced studies by utilizing well-defined mathematical properties of such polynomials, such as linear transformations, interpolations, continuities, differentiabilities, etc. Next, we consider a \textit{bivariate polynomial problem for matrices}, stated as follows:

$\bullet$ \textit{For a given matrix $({a_{ij}}) \in \mathcal{R}^{m \times n}$ and an $mn$-dimensional subspace $\mathcal{P} \subset \Pi^2$, find a polynomial $p \in \mathcal{P}$, such that 
\begin{equation}
    p(i,j) = {a_{ij}} \textrm{ for all } (i,j) \in \mathcal{X}=\{1,2,\ldots,m\}\times \{1,2,\ldots,n\} \label{DPPM}
\end{equation}}
is named as Dharm polynomial problem for matrices (DPPM). As an advanced application, the DPPM \eqref{DPPM} provides a \textit{sufficient condition} for the map $\mathcal{D}_{p}:\mathcal{R}^{m \times n} \rightarrow \mathcal{P}$, defined as 
\begin{equation}\label{dp isomorphism}
    \mathcal{D}_{p}{(A)} = p_{A}(x,y) \;\textrm{for all}\; A \in \mathcal{R}^{m \times n},
\end{equation}
to be an \textit{isomorphism} \cite{9}, where $p_A \in \mathcal{P}$ denoes the unique polynomial satisfying the DPPM \eqref{DPPM}. 

This article focuses on the existence, uniqueness, and construction of the polynomials satisfying the DPPM \eqref{DPPM} in some $mn$-dimensional subspace of $\Pi^2$ of various total degrees and shows the possible applicability of results in bivariate polynomial interpolation problems concerning the rectangular schemes of bivariate interpolation points on the natural Cartesian grid $\mathcal{X}$.

The rest of the article is organized as follows. Section 2 deals with the basic concepts and preliminaries of the article. Section 3 focuses on the results of the article. Section 4 discusses the possible applicability of the DPPM \eqref{DPPM} in bivariate polynomial interpolation problems. Section 5 considers some numerical examples to justify the theoretical findings. At last, section 6 concludes the article.

\section{Basic Concepts and Preliminaries}\label{sec2}
For a given vector $\mathcal{F}= \left(f({x_i},y_j):i=1,2,\ldots,m; j=1,2,\ldots,n\right)$ of prescribed real values and a subspace $\mathcal{Q} \subset \Pi^2$, find a polynomial $p \in \mathcal{Q}$, for a finite set $\Theta = \left\{\left(x_i,y_j\right):i=1,2,\ldots,m; j=1,2,\ldots,n \right\}$ of $mn$ pair-wise distinct interpolation points in $\mathcal{R}^{2}$, such that
\begin{equation}\label{lag:int}
    p(\Theta) = \mathcal{F}, \;i.e., \; p(x_i,y_j) = f(x_{i},y_j) \textrm{ for all }\; (x_i,y_j) \in \Theta
\end{equation}
represents the \textit{Lagrange's bivariate polynomial interpolation problem}. In this case, the interpolation points $\left(x_i,y_j\right)$, $(x_i,y_j) \in \Theta$ are also called \textbf{\textit{nodes}}. The LBVPIP \eqref{lag:int} for $\Theta$ is called \textbf{\textit{correct}} in $\mathcal{Q}$, if there exists a unique polynomial $p \in \mathcal{Q}$ such that \eqref{lag:int} holds for any given vector $\mathcal{F}$. For this, it is necessary that $\dim(\mathcal{Q})=mn$. Also, if the LBVPIP \eqref{lag:int} with respect to $\Theta$ is correct in $\mathcal{Q}$, then $\mathcal{Q}$ is called \textbf{\textit{correct space}} or \textbf{\textit{interpolation space}} for the LBVPIP \eqref{lag:int} with respect to $\Theta$, cf. \cite{10,11}. 

Without loss of generality, for $(a_{ij}) = \left(f(x_i,y_j)\right) \in \mathcal{R}^{m \times n}$ and the subspace $\mathcal{P} \subset \Pi^2$, the LBVPIP \eqref{lag:int} with respect to $\mathcal{X}$ is equivalent to the DPPM \eqref{DPPM}. Thus, if $\mathcal{P}$ is a correct space for the LBVPIP \eqref{lag:int} for $\mathcal{X}$, then there must exist a unique polynomial $p_A \in \mathcal{P}$ that satisfies the DPPM \eqref{DPPM} for any given choice of the matrix $A \in \mathcal{R}^{m \times n}$ and the map \eqref{dp isomorphism} can always be defined. In other words, there exist a unique polynomial $p_A \in \mathcal{P}$ that satisfies the DPPM \eqref{DPPM} for any choice of the matrix $A \in \mathcal{R}^{m \times n}$, if the matrix 
\begin{equation}\label{sample matrix}
    \begin{pmatrix} u(i,j): (i,j) \in \mathcal{X}, u \in \mathcal{B} \end{pmatrix} \in \mathcal{R}^{mn\times mn}
\end{equation}
is non-singular for a basis $\mathcal{B}$ of $\mathcal{P}$. Thus, the DPPM \eqref{DPPM}, equivalently the LBVPIP \eqref{lag:int} for $\mathcal{X}$ is correct in $\mathcal{P}$ if the determinant of the matrix \eqref{sample matrix} is non-zero for any choice of the basis $\mathcal{B}$ of $\mathcal{P}$. The interpolation theory in \cite{10,11,12} tells that any rectangular scheme of $m$ by $n$ bivariate nodes on the natural Cartesian grid $\mathcal{X}$ is $\Pi_{k}^{2}$-independent if and only if $k=m+n-2$, but then it cannot be $\Pi_{k}^{2}$-complete. Therefore, the DPPM \eqref{DPPM} cannot be $\Pi_{k}^{2}$-correct for any $k$ (the only exception being the trivial case $m = n = 1$: single node is of course 0-correct). 

In \cite{13}, the author Narumi correctly claims that the LBVPIP \eqref{lag:int} with respect to $\mathcal{X}$ can always be correct in a standard $mn$-dimensional tensor product subspace $\mathcal{P}_{m}^{n}=\Pi_{m-1} \otimes \Pi_{n-1} \subset \Pi^2_{m+n-2}$ \cite{14}, such that, if $p \in \mathcal{P}_{m}^{n}$, then
\begin{equation}\label{dpsosk} 
p(x,y)= {\sum_{k_1=0}^{m-1} {\sum_{k_2=0}^{n-1}}}{\lambda _{{k_1,}{k_2}}}{x^{k_1}}{y^{k_2}}, \;  {\lambda _{{k_1,}{k_2}}} \in  \mathcal{R} \; \textrm{for all} \; 0 \leq k_1 \leq {m-1}, 0 \leq k_2 \leq {n-1}, 
\end{equation}  
where $\left\{x^{k_1}y^{k_2}:0 \leq k_1 \leq m-1, 0 \leq k_2 \leq n-1 \right\}$ is the set of standard basis of $\mathcal{P}_{m}^{n}$. Moreover, a formula has been submitted to construct the coefficients of the required polynomial, involving Newton-divided differences of the data values. However, the proof of the existence and uniqueness of the needed polynomial in $\mathcal{P}_{m}^{n}$ has not been included. The review article \cite{12} guessed that the author Narumi might have pursued the tensor product approach as that time, the tensor product of univariate interpolating polynomials when the bivariate interpolation points lie on a natural Cartesian product grid $\mathcal{X}$, was the classical approach to bivariate polynomial interpolation, see also \cite{15}. In general, $\mathcal{P}_{m}^{n}$ is used as a standard interpolation space for the LBVPIP \eqref{lag:int} for $\mathcal{X}$. 

Moreover, given the Kergin interpolation \cite{16}, there is at least one $mn$-dimensional subspace of $\Pi_{mn-1}^2$ in which the LBVPIP \eqref{lag:int} for $\mathcal{X}$ can always be correct. In general, there may possibly exist infinitely many linearly independent $mn$-dimensional subspaces of $\Pi^{2}$, of various total degrees, in which the LBVPIP \eqref{lag:int} for $\mathcal{X}$ can always be correct. Indeed, the identification or construction of such correct spaces is much more difficult in case of more than one variable, see \cite{1,10,11,12} for details. Consequently, the DPPM \eqref{DPPM} can always be correct in $\mathcal{P}_{m}^{n}$ and there may exist one or more linearly independent $mn$-dimensional subspaces of $\Pi^2$, of various total degrees, in which the DPPM \eqref{DPPM} can always be correct. 

\section{Main Results}\label{sec3}
\begin{thm}\label{sufficient condition}
If $\mathcal{P}$ is a correct space for the DPPM \eqref{DPPM}, then the map $\mathcal{D}_{p}:\mathcal{R}^{m \times n} \rightarrow \mathcal{P}$, defined as \eqref{dp isomorphism}, is an isomorphism.
\end{thm}
\begin{proof}
Let $\mathcal{P}$ be a correct space for the DPPM \eqref{DPPM}, then $\dim(\mathcal{P})=mn$ and there always exist a unique polynomial $p_{A} \in \mathcal{P}$ that satisfies the DPPM \eqref{DPPM} for each $A=({a_{ij}})\in \mathcal{R}^{m \times n}$, i.e., the map \eqref{dp isomorphism} is well-defined. Now, let ${B}=({b_{ij}}) \in \mathcal{R}^{m \times n}$, ${C} =(c_{ij}) \in \mathcal{R}^{m \times n}$, and $\lambda$ be a real scalar. Then, ${B} + \lambda {C} = (b_{ij}+ \lambda c_{ij}) \in \mathcal{R}^{m \times n}$. Therefore, there exists $p_{B} \in \mathcal{P}$, $p_{C} \in \mathcal{P}$, and $p_{B+\lambda C} \in \mathcal{P}$, such that 
\begin{eqnarray}
&& p_{B}(i,j) = b_{ij} \; \textrm{for all} \; (i,j) \in \mathcal{X}, \label{PA:2.1}\\
&& p_{C}(i,j) = c_{ij} \; \textrm{for all} \; (i,j) \in \mathcal{X}, \label{PB:2.1}\\
&& \textrm{and } p_{B+\lambda C}(i,j) = b_{ij} + \lambda c_{ij} \; \textrm{for all} \; (i,j) \in \mathcal{X}. \label{DPPM:addition2.1} 
\end{eqnarray}
From \eqref{PA:2.1}, \eqref{PB:2.1} and \eqref{DPPM:addition2.1}, we have
\begin{eqnarray*}
&& p_{B+\lambda C}(i,j) = p_{B}(i,j)+\lambda p_{C}(i,j) \; \textrm{for all} \; (i,j) \in \mathcal{X}, \label{PA:2.12}
\end{eqnarray*}
i.e., the map \eqref{dp isomorphism} is linear. Again, for each $p \in \mathcal{P}$, there exists $(p(i,j)) \in \mathcal{R}^{m \times n}$, i.e., the map \eqref{dp isomorphism} is surjective. Since, $\dim(\mathcal{R}^{m \times n})=\dim(\mathcal{P})$, thus the map \eqref{dp isomorphism} is bijective and consequently invertible. Hence, the proof is completed.
\end{proof}

\begin{rem}\label{inverse: isomorphic map def}
The inverse map $\mathcal{D}_{p}^{-1}:\mathcal{P}\rightarrow \mathcal{R}^{m \times n}$, of the map \eqref{dp isomorphism}, is given as
\begin{equation}\label{inverse isomorphic map2.1}
     \mathcal{D}_{p}^{-1}\left({p(x,y)}\right) = {\left(p(i,j)\right)}\in \mathcal{R}^{m \times n} \; \textrm{for all} \; p \in \mathcal{P}.
\end{equation}
\end{rem}

A straightforward consequence of Theorem \ref{sufficient condition} is as follows.
\begin{cor}\label{special property}
If $\mathcal{P}$ is a correct space for the DPPM \eqref{DPPM} and $p, q \in \mathcal{P}$, then 
\begin{enumerate}
    \item $p(x,y)=0$ in $\mathcal{R}^2$ if and only if $p(x,y)=0$ in $\mathcal{X}$, and
    \item $p(x,y)=q(x,y)$ in $\mathcal{R}^2$ if and only if $p(x,y)=q(x,y)$ in $\mathcal{X}$.
\end{enumerate}
\end{cor}

\begin{rem}\label{remark2.1}
If $\mathcal{P}$ is a correct space for the DPPM 
 \eqref{DPPM}, then the number of roots of a polynomial $p \in \mathcal{P}$ can never exceed $mn$ and two polynomials $p,q \in \mathcal{P}$ are said to be equal if and only if $p(x,y)=q(x,y)$ in $\mathcal{X}$.
\end{rem}

\begin{thm} \label{thm:dp} 
For every $(a_{ij}) \in \mathcal{R}^{m \times n}$, there exists a unique $p \in \mathcal{P}_{m}^{n}$ that satisfies the DPPM \eqref{DPPM}.
\end{thm}

\begin{proof}
The proof consists two parts, existence and uniqueness. To prove this, it is sufficient to show that the matrix \eqref{sample matrix} is non-singular for the basis $\mathcal{B} = \left\{x^{k_1}y^{k_2}:0 \leq k_1 \leq m-1, 0 \leq k_2 \leq n-1 \right\}$. For better understanding, we will prove existence and uniqueness independently.

\noindent \textit{Existence:} Let $(a_{ij}) \in \mathcal{R}^{m \times n}$ is a given matrix. For the set of interpolation points 
	\begin{equation}\label{rthrow:interpolation points}
            \Omega_n = \left\{j:j=1,2,\ldots,n\right\} 
	\end{equation}
concerning the $k$th row of the given matrix $(a_{ij}) \in \mathcal{R}^{m \times n}$, there exists a univariate polynomial $p_k \in \Pi_{n-1}$ that satisfy 
	\begin{equation}\label{intp:1}
            p_k(j)=a_{kj}\textrm{ for all } j=1,2,\ldots,n
	\end{equation}
for all $k=1,2,\ldots,m$. Thus, $(a_{ij}) \in \mathcal{R}^{m \times n}$ can be represented as $(a_{ij})=\left(p_k \in \Pi_{n-1}: k=1,2,\ldots,m\right)_{m \times 1}$. Again, for the set of interpolation points 
	\begin{equation}\label{gen:interpolation points}
            \Omega_m = \left\{i:i=1,2,\ldots,m\right\}  
	\end{equation}
there exists a bivariate polynomial $p \in \Pi_{m-1} \times \Pi_{n-1}$ that satisfy 
	\begin{equation}\label{intp:2}
            p(i,y)=p_i(y) \textrm{ for all } i=1,2,\ldots,m.
	\end{equation}
Consequently, there exist a polynomial $p \in \mathcal{P}_{m}^{n}$ that satisfy the DPPM \eqref{DPPM} for any given choice of the matrix $(a_{ij}) \in \mathcal{R}^{m \times n}$. 

\noindent \textit{Uniqueness:} Let $p \in \mathcal{P}_{m}^{n}$ be the polynomial, of the form \eqref{dpsosk}, that satisfies the DPPM \eqref{DPPM}. Then,
\begin{equation} \label{uniqeness:2.1}
	{\sum_{k_1=0}^{m-1} {\sum_{k_2=0}^{n-1}}}{\lambda _{{k_1,}{k_2}}}{i^{k_1}}{j^{k_2}}= {a_{ij}} \; \textrm{for all} \; (i,j) \in \mathcal{X}.
\end{equation}
The system of equations \eqref{uniqeness:2.1} can be written as
\begin{equation}\label{sys:gen2.1}
    {\Xi}\mathcal{G}=\nu,
\end{equation}
where $\nu = \left( \begin{array}{ccccccccccc}
a_{11} & a_{12} & \ldots & a_{1n} & \ldots & \ldots & a_{m1} & a_{m2} & \ldots & a_{mn}
\end{array} \right) \in \mathcal{R}^{mn \times 1}$,
\begin{flalign*}
\mathcal{G}=\left( \begin{array}{ccccccccccc}
\lambda_{0,0} & \lambda_{0,1} & \ldots & \lambda_{0,n-1} & \ldots & \ldots & \lambda_{m-1,0} & \lambda_{m-1,1} & \ldots & \lambda_{m-1,n-1}
\end{array} \right) \in \mathcal{R}^{mn \times 1}, &&
\end{flalign*}
and
\begin{flalign*}
& \Xi= \left( \begin{array}{cccccccccc}
{1^{0}}{1^{0}} & {1^{0}}{1^{1}} & \ldots & {1^{0}}{1^{n-1}} & \ldots & {1^{m-1}}{1^0} & {1^{m-1}}{1^{1}} & \ldots & {1^{m-1}}{1^{n-1}} \\
{1^{0}}{2^{0}} & {1^{0}}{2^{1}} & \ldots & {1^{0}}{2^{n-1}} & \ldots & 1^{m-1}{2^{0}} & {1^{m-1}}{2^{1}} & \ldots & {1^{m-1}}{2^{n-1}} \\
\vdots & \vdots & \ldots & \vdots & \ldots & \vdots & \vdots & \ldots & \vdots\\
{1^{0}}{n^{0}} & {1^{0}}{n^{1}} & \ldots & {1^{0}}{n^{n-1}} & \ldots & {1^{m-1}}{n^0} & {1^{m-1}}{n^{1}} & \ldots & {1^{m-1}}{n^{n-1}} \\
{2^{0}}{1^{0}} & {2^{0}}{1^{1}} & \ldots & {2^{0}}{1^{n-1}} & \ldots & 2^{m-1}{1^{0}} & {2^{m-1}}{1^{1}} & \ldots & {2^{m-1}}{1^{n-1}} \\
{2^{0}}{2^{0}} & {2^{0}}{2^{1}} & \ldots & {2^{0}}{2^{n-1}} & \ldots & 2^{m-1}{2^{0}} & {2^{m-1}}{2^{1}} & \ldots & {2^{m-1}}{2^{n-1}} \\
\vdots & \vdots & \ldots & \vdots & \ldots & \vdots & \vdots & \ldots & \vdots\\
{2^{0}}{n^{0}} & {2^{0}}{n^{1}} & \ldots & {2^{0}}{n^{n-1}} & \ldots & {2^{m-1}}{n^0} & {2^{m-1}}{n^{1}} & \ldots & {2^{m-1}}{n^{n-1}} \\
\vdots & \vdots & \ldots & \vdots & \ldots & \vdots & \vdots & \ldots & \vdots\\
{m^{0}}{1^{0}} & {m^{0}}{1^{1}} & \ldots & {m^{0}}{1^{n-1}} & \ldots & m^{m-1}{1^{0}} & {m^{m-1}}{1^{1}} & \ldots & {m^{m-1}}{1^{n-1}} \\
{m^{0}}{2^{0}} & {m^{0}}{2^{1}} & \ldots & {m^{0}}{2^{n-1}} & \ldots & m^{m-1}{2^{0}} & {m^{m-1}}{2^{1}} & \ldots & {m^{m-1}}{2^{n-1}}  \\
\vdots & \vdots & \ldots & \vdots & \ldots & \vdots & \vdots & \ldots & \vdots\\
{m^{0}}{n^{0}} & {m^{0}}{n^{1}} & \ldots & {m^{0}}{n^{n-1}} & \ldots & m^{m-1}{n^{0}} & {m^{m-1}}{n^{1}} & \ldots & {m^{m-1}}{n^{n-1}} 
\end{array} \right)_{{mn} \times {mn}}.
\end{flalign*}
Using the Kronecker (or tensor) product of matrices \cite{17}, one may get
\begin{flalign*}
\Xi & = \left( \begin{array}{cccc}
{1^{0}} & {1^{1}} & \ldots & {1^{m-1}} \\
{2^{0}} & {2^{1}} & \ldots & {2^{m-1}} \\
\vdots & \vdots & \ldots & \vdots \\
{m^{0}} & {m^{1}} & \ldots & {m^{m-1}} 
\end{array} \right)_{m\times m} \otimes \left( \begin{array}{cccc}
{1^{0}} & {1^{1}} & \ldots & {1^{n-1}} \\
{2^{0}} & {2^{1}} & \ldots & {2^{n-1}} \\
\vdots & \vdots & \ldots & \vdots \\
{n^{0}} & {n^{1}} & \ldots & {n^{n-1}} 
\end{array} \right)_{n\times n} \\
\textrm{and } \det(\Xi) & =  \left(\begin{vmatrix}
{1^{0}} & {1^{1}} & \ldots & {1^{m-1}} \\
{2^{0}} & {2^{1}} & \ldots & {2^{m-1}} \\
\vdots & \vdots & \ldots & \vdots \\
{m^{0}} & {m^{1}} & \ldots & {m^{m-1}} 
\end{vmatrix}\right)^{m} \times \left(\begin{vmatrix}
{1^{0}} & {1^{1}} & \ldots & {1^{n-1}} \\
{2^{0}} & {2^{1}} & \ldots & {2^{n-1}} \\
\vdots & \vdots & \ldots & \vdots \\
{n^{0}} & {n^{1}} & \ldots & {n^{n-1}} 
\end{vmatrix}\right)^{n}. &&
\end{flalign*}
Therefore, using the standard Vandermonde determinant \cite{18}, we get $\det({\Xi}) \neq 0$, i.e., the system of equations \eqref{sys:gen2.1} possesses a unique solution. Hence, the proof is completed.
\end{proof}

\begin{rem}
Using the system of equations \eqref{sys:gen2.1}, the coefficients of the required unique polynomial $p \in \mathcal{P}_{m}^{n}$ that satisfies the DPPM \eqref{DPPM} for each $(a_{ij}) \in \mathcal{R}^{m \times n}$ can be determined by setting ${\mathcal{G}}={\Xi^{-1}}{\nu}$.
\end{rem}

The proof of existence part of Theorem \ref{thm:dp} provides an algorithm to construct the polynomial $p \in \mathcal{P}_{m}^{n}$ that satisfies the DPPM \eqref{DPPM} for each $(a_{ij}) \in \mathcal{R}^{m \times n}$, considering univariate polynomial interpolation formulas \cite{19,20}. Next, we present the following results for the algorithm given in the existence part of Theorem \ref{thm:dp}.

\begin{cor} \label{cor:construction by lagrange} 
For every $(a_{ij}) \in \mathcal{R}^{m \times n}$, there exists a unique $p \in \mathcal{P}_{m}^{n}$ that satisfies the DPPM \eqref{DPPM} is given as
	\begin{equation} \label{lag:mat}
		p(x,y)=\sum_{k=1}^{m} \left(\prod_{\alpha=1, \alpha \neq k}^{m} \left({\frac {x-\alpha}{k-\alpha}}\right) \right){p_k(y)}, 
	\end{equation}
where
	\begin{equation} \label{lag:row}
        \quad p_k(y)=\sum_{r=1}^{n} \left(\prod_{\alpha=1, \alpha \neq r}^{n} \left({\frac {y-\alpha}{r-\alpha}}\right) \right){a_{kr}} \textrm{ for all } k=1,2,\ldots,m.
	\end{equation}
\end{cor}

\begin{proof}
Using univariate Newton-Lagrange's polynomial interpolation formula, the polynomials $p_k \in \Pi_{n-1}$, $k=1,2,\ldots,m$ for the set of interpolation points \eqref{rthrow:interpolation points} that satisfy the problem \eqref{intp:1} can be written as 
\begin{equation} \label{lag:1}
        p_k(y)=l_{k,1}(y)a_{k1}+l_{k,2}(y)a_{k2}+\ldots+l_{k,n}(y)a_{kn},
\end{equation}
where $l_{k,r}(y)$, $r=1,2,\ldots,n$ are the univariate polynomials in $\Pi_{n-1}$. The polynomials \eqref{lag:1} satisfy \eqref{intp:1}, if and only if 
\begin{equation} \label{lag:2}
        l_{k,r}(y_j)=\begin{cases}
                    0,  & \text{if } r \neq j\\
                    1,  & \text{if } r = j 
                \end{cases} \textrm{ for all }r=1,2,\ldots,n.
\end{equation}
Therefore, the polynomials $l_{k,r}(y)$, $r=1,2,\ldots,n$ that satisfy \eqref{lag:2} can be written as
\begin{equation} \label{lag:3}
        l_{k,r}(y)= \prod_{\alpha=1, \alpha \neq r}^{n} {\frac {(y-\alpha)}{(r-\alpha)}} \textrm{ for all } r=1,2,\ldots,n.
\end{equation}
Thus, using \eqref{lag:3} in \eqref{lag:1}, we get \eqref{lag:row}. Again, using univariate Newton-Lagrange's polynomial interpolation formula, the polynomial $p \in \Pi_{m-1} \times \Pi_{n-1}$ for the set interpolation points \eqref{gen:interpolation points} that satisfy the problem \eqref{intp:2} can be written as 
\begin{equation} \label{lag:4}
        p(x,y)=p_1(y)L_1(x)+p_2(y)L_2(x)+\ldots+p_m(y)L_m(x),
\end{equation}
where $L_k(x)$, $k=1,2,\ldots,m$ are the univariate polynomials in $\Pi_{m-1}$. The polynomial \eqref{lag:4} satisfy \eqref{intp:2}, if and only if 
\begin{equation} \label{lag:5}
        L_k(x_i)=\begin{cases}
                    0,  & \text{if } k \neq i\\
                    1,  & \text{if } k = i 
                \end{cases} \textrm{ for all } k=1,2,\ldots,m.
\end{equation}
and the polynomials $L_k(x)$ that satisfy the conditions \eqref{lag:5} can be written as
\begin{equation} \label{lag:6}
        L_k(x)= \prod_{\alpha = 1, \alpha \neq k}^{m} {\frac {(x-\alpha)}{(k-\alpha)}} \textrm{ for all } k=1,2,\ldots,m.
\end{equation}
Therefore, combining \eqref{lag:4} and \eqref{lag:6}, we get \eqref{lag:mat}. Clearly, $p \in \mathcal{P}_{m}^{n}$ and the uniqueness followed from Theorem \ref{thm:dp}. This completes the result.
\end{proof}

\begin{cor} \label{cor:construction by newton forward} 
For every $(a_{ij}) \in \mathcal{R}^{m \times n}$, there exists a unique $p \in \mathcal{P}_{m}^{n}$ that satisfies the DPPM \eqref{DPPM} is given as
\begin{equation}\label{gnf:3}
	p(x,y)= \sum_{r=0}^{m-1}{{{x}\choose{r}}\Delta ^r {p_1(y)}},
\end{equation}
where
\begin{equation} \label{gnf:4}
    \Delta^{r} p_1(y)=\sum_{\beta=0}^{r}{(-1)^{\beta}}{r \choose \beta}{p_{r+1-\beta}(y)} \textrm{ for all } r=1,2,\ldots,m-1
\end{equation}
\newline
and
\begin{equation} \label{gnf:1}
\begin{aligned}
        p_k(y) = {} \sum_{r=0}^{n-1}{{{y}\choose{r}}\Delta ^r {a_{k1}}} \textrm{ for all } k=1,2,\ldots,m,
\end{aligned}
\end{equation}
where
\begin{equation} \label{gnf:2}
    \Delta^{r} a_{k1}=\sum_{\alpha=0}^{r}{(-1)^{\alpha}}{r \choose \alpha}{a_{k(r+1-\alpha)}} \textrm{ for all } r=1,2,\ldots,n-1.
\end{equation}
\end{cor}

\begin{proof}
Using univariate Newton-forward difference polynomial interpolation formula instead of univariate Newton-Lagrange's polynomial interpolation formula in the proof of Corollary \ref{cor:construction by lagrange} completes the result.
\end{proof}

Next, a class of $mn$-dimensional subspaces of $\Pi^2_{s(mn-1)}$ is considered, involving two real parameters $\alpha$ and $\beta$, and it is proven that the DPPM \eqref{DPPM} can be correct in these subspaces for uncountably many choices of the pair $(\alpha,\beta)$. 
\begin{defn}\label{new:subspace}
For some real scalars $\alpha$, $\beta$ and $s \in \mathbb{N}$, ${}_{\alpha}^{\beta}\Pi_{s(mn-1)}^{2}$ is the space of bivariate polynomials with real coefficients, of total degree up to $s(mn-1)$, such that if $p \in {_{\alpha}^{\beta}\Pi_{s(mn-1)}^{2}}$, then 
\begin{eqnarray}\label{polynomial form}
       p(x,y)= \sum_{k=0}^{mn-1}{{\lambda_{k}}(\alpha x^s + \beta y^s)^{k}}, \quad  {\lambda _{k}} \in  \mathcal{R} \;\; \textrm{for all} \;\; 0 \leq k \leq {mn-1}.
\end{eqnarray}
where $\lambda _{k}$ are some real scalars for all $0 \leq k \leq {mn-1}$.
\end{defn}

\begin{rem}\label{basis2}
The space ${}_{\alpha}^{\beta}\Pi_{s(mn-1)}^{2}$ is an $mn$-dimensional subspace of $\Pi_{s(mn-1)}^2$ as $\left({\mu}u+v\right) \in {}_{\alpha}^{\beta}\Pi_{s(mn-1)}^{2}$ for all $u,v \in {}_{\alpha}^{\beta}\Pi_{s(mn-1)}^{2}$ and scalar $\mu$, where a set of standard bases of ${}_{\alpha}^{\beta}\Pi_{s(mn-1)}^{2}$ is $\left\{(\alpha x^s + \beta y^s)^{k}:0 \leq k \leq {mn-1}\right\}$.
\end{rem}

\begin{thm}\label{thm:gen1.1} 
For every $({a_{ij}}) \in \mathcal{R}^{m \times n}$, there exists a unique $p \in {_{\alpha}^{\beta}\Pi_{s(mn-1)}^{2}}$ that satisfies the DPPM \eqref{DPPM}, provided the pair $(\alpha, \beta)$ satisfy the condition 
\begin{equation}\label{cond:nodes}
    \alpha i_{1}^{s} + \beta j_{1}^{s} \neq \alpha i_{2}^{s} + \beta j_{2}^{s} \textrm{ for all } (i_1,j_1) \neq (i_2,j_2)
\end{equation}
where $(i_1,j_1),(i_2,j_2) \in \mathcal{X}$. 
\end{thm}

\begin{proof}
Consider the transformation $\phi_s:\mathcal{X} \rightarrow \mathcal{R}$, defined as
\begin{equation} \label{pf:11}
		\phi_s(i,j)=\alpha i^s + \beta j^s \; \textrm{for all} \; (i,j) \in \mathcal{X}.
	\end{equation}
If the pair $(\alpha, \beta)$ satisfies the condition \eqref{cond:nodes}, the transformation \eqref{pf:11} is a injective that projects all of bivariate nodes of $\mathcal{X}$ onto a set of $mn$-univariate nodes to the curve orthogonal to $\alpha x^s + \beta y^s = 0$ in such a manner that all projections are different. Thus, using \eqref{pf:11}, the DPPM \eqref{DPPM} can be transformed into a univariate polynomial interpolation problem of the form
\begin{eqnarray}\label{DPPM:modified1}
            p(i,j) = \widetilde{p}(\phi_s(i,j)) = \widetilde{p}(\alpha i^s + \beta j^s) = {a_{ij}} \textrm{ for all } (i,j) \in \mathcal{X}.
\end{eqnarray}
Consequently, using a univariate polynomial interpolation formula for the interpolation problem \eqref{DPPM:modified1} for $\mathcal{X}$, there exists a polynomial $p \in {} _{\alpha}^{\beta}\Pi_{s(mn-1)}^2$ that satisfy the DPPM \eqref{DPPM}. Again, let $p \in {_{\alpha}^{\beta}\Pi_{s(mn-1)}^{2}}$ be the polynomial, as defined in \eqref{polynomial form}, that satisfies the DPPM \eqref{DPPM} for the given $({a_{ij}}) \in \mathcal{R}^{m \times n}$. Then
\begin{equation} \label{uni:gen22}
	\sum_{k=0}^{mn-1}{{\lambda_{k}}(\alpha i^s + \beta j^s)^{k}} = {a_{ij}} \; \textrm{for all} \; (i,j) \in \mathcal{X},
\end{equation}
i.e., the coefficients must satisfy the system of equations of the form
\begin{equation}\label{sys:gen1}
    {\Lambda}\mathcal{K}=\mu,
\end{equation}
where, $\mathcal{K} = 
	{\begin{pmatrix}			
		\lambda_0 ~
		\lambda_1~
		\ldots~
		\lambda_{n-1}~
		\lambda_n~
		\ldots~
		\lambda_{mn-1}
		\end{pmatrix}} \in \mathcal{R}^{mn \times 1}$, 
$\mu =
	{\begin{pmatrix}			
	    a_{11}~
		a_{12}~
		\ldots~
		a_{1n}~
		a_{21}~
		\ldots~
		a_{mn}
		\end{pmatrix}} \in \mathcal{R}^{mn \times 1}$ and
\begin{flalign*} \label{coefficient:matrix2}
    	{\Lambda} = {\begin{pmatrix}			
		1 & (\alpha + \beta) & (\alpha + \beta)^2 & \ldots & (\alpha + \beta)^{mn-1}\\
		1 & (\alpha + 2^s\beta) & (\alpha + 2^s\beta)^2 & \ldots & (\alpha + 2^s\beta)^{mn-1}\\
		\vdots & \vdots & \vdots &  & \vdots\\
		1 & (\alpha + n^s\beta) & (\alpha + n^s\beta)^2 & \ldots & (\alpha + n^s\beta)^{mn-1}\\
		1 & (2^s\alpha + \beta) & (2^s\alpha + \beta)^2 & \ldots & (2^s\alpha + \beta)^{mn-1}\\
		\vdots & \vdots & \vdots &  & \vdots\\
		1 & (m^s\alpha + n^s\beta) & (m^s\alpha + n^s\beta)^2 & \ldots & (m^s\alpha + n^s\beta)^{mn-1}\\
		\end{pmatrix}}_{mn \times mn}.
\end{flalign*}
Therefore, if the pair $(\alpha, \beta)$ satisfies the condition \eqref{cond:nodes}, then $\det(\Lambda) \neq 0$, using Vandermonde determinant. Thus, the system of equations \eqref{sys:gen1} has a unique solution. Hence, the proof is completed.
\end{proof}

\begin{rem}
The condition \eqref{cond:nodes} holds for uncountably many choices of the pair $(\alpha, \beta)$ for every $s \in \mathbb{N}$. This establishes the existence of a class of $mn$-dimensional subspaces of $\Pi^2$, of the form ${}_{\alpha}^{\beta}\Pi_{s(mn-1)}^{2}$, for each $s \in \mathbb{N}$, in which the DPPM \eqref{DPPM} can always be correct. For a positive integer $s$, the coefficients of the required polynomial $p \in {}_{\alpha}^{\beta}\Pi_{s(mn-1)}^{2}$ that uniquely satisfies the DPPM \eqref{DPPM} can also be determined from the system of equation \eqref{sys:gen1} using $\mathcal{K} = {\Lambda}^{-1}{\mu}$, provided the pair $(\alpha, \beta)$ satisfies the condition \eqref{cond:nodes}. 
\end{rem}

\begin{rem}
The system \eqref{uni:gen22} is not the only one that gives the solution to the DPPM \eqref{DPPM}. There is another:
\begin{equation}\label{sys:gen1r}
    \bar{\Lambda}{\mathcal{K}}=\bar{\mu},
\end{equation}
where $\bar{\mu} =
	{\begin{pmatrix}			
	    a_{11}~
		a_{21}~
		\ldots~
		a_{m1}~
		a_{12}~
		\ldots~
		a_{mn}
		\end{pmatrix}} \in \mathcal{R}^{mn \times 1}$ and
\begin{flalign*} \label{coefficient:matrix2r}
    	\bar{\Lambda} = {\begin{pmatrix}			
		1 & (\alpha + \beta) & (\alpha + \beta)^2 & \ldots & (\alpha + \beta)^{mn-1}\\
		1 & (2^s\alpha + \beta) & (2^s\alpha + \beta)^2 & \ldots & (2^s\alpha + \beta)^{mn-1}\\
		\vdots & \vdots & \vdots &  & \vdots\\
		1 & (m^s\alpha + \beta) & (m^s\alpha + \beta)^2 & \ldots & (m^s\alpha + \beta)^{mn-1}\\
		1 & (\alpha + 2^s\beta) & (\alpha + 2^s\beta)^2 & \ldots & (\alpha + 2^s\beta)^{mn-1}\\
		\vdots & \vdots & \vdots &  & \vdots\\
		1 & (m^s\alpha + n^s\beta) & (m^s\alpha + n^s\beta)^2 & \ldots & (m^s\alpha + n^s\beta)^{mn-1}\\
		\end{pmatrix}}_{mn \times mn}.
\end{flalign*}
\end{rem}

The proof of Theorem \ref{thm:gen1.1} also provides an algorithm to construct the polynomial $p \in {}_{\alpha}^{\beta}\Pi_{mn-1}^{2}$ that uniquely satisfies the DPPM \eqref{DPPM}, provided the pair $(\alpha, \beta)$ satisfies the condition \eqref{cond:nodes} for each $s \in \mathbb{N}$. On applying the univariate Newton-divided difference polynomial interpolation formula for the interpolation problem \eqref{DPPM:modified1} for $\mathcal{X}$, we get the following result.

\begin{cor}\label{corollary:1.0}
For every $({a_{ij}}) \in \mathcal{R}^{m \times n}$ and positive integer $s$, there exists a unique $p \in {_{\alpha}^{\beta}\Pi_{s(mn-1)}^{2}}$ that satisfies the DPPM \eqref{DPPM}, such that the pair $(\alpha, \beta)$ satisfy the condition \eqref{cond:nodes}, is given as 
\begin{equation}\label{pol:common}
\begin{aligned}
	p(x,y)= & {} \widetilde{p}(\phi_s(x,y)) = \widetilde{p}(\alpha x^s + \beta y^s) \\ = & \;  \widetilde{a}_{\alpha+\beta} + {((\alpha x^s + \beta y^s)-(\alpha + \beta))}{F[\widetilde{a}_{\alpha+\beta},\widetilde{a}_{\alpha+\beta 2^s}]} \\
	& \; \; +{((\alpha x^s + \beta y^s)-(\alpha + \beta))}{((\alpha x^s + \beta y^s)-(\alpha + \beta 2^s))}{F[\widetilde{a}_{\alpha+\beta},\widetilde{a}_{\alpha+\beta 2^s},\widetilde{a}_{\alpha+\beta 3^s}]}   \\ 
	& \qquad + ... + {((\alpha x^s + \beta y^s)-(\alpha + \beta))}{((\alpha x^s + \beta y^s)-(\alpha + \beta 2^s))}\ldots\\
	& \qquad \qquad {((\alpha x^s + \beta y^s)-(\alpha m^s+ \beta (n-1)^s))}{F[\widetilde{a}_{\alpha+\beta},\widetilde{a}_{\alpha+\beta 2^s},\ldots,\widetilde{a}_{\alpha m^s+\beta n^s}]},
\end{aligned}
\end{equation}
where  
\begin{equation*}
    {F[\widetilde{a}_{\alpha+\beta},\widetilde{a}_{\alpha+\beta 2^s}]} = \frac{{F[\widetilde{a}_{\alpha+\beta 2^s}]}-{F[\widetilde{a}_{\alpha+\beta}]}}{({\alpha +\beta 2^s})-({\alpha +\beta})}
\end{equation*}
and \\ ${F[\widetilde{a}_{\alpha+\beta},\widetilde{a}_{\alpha+\beta 2^s},\ldots,\widetilde{a}_{\alpha m^s+\beta n^s}]}$
\begin{equation*}
   = \frac{{F[\widetilde{a}_{\alpha+\beta 2^s},\widetilde{a}_{\alpha+\beta 3^s},\ldots,\widetilde{a}_{\alpha m^s+\beta n^s}]} - {F[\widetilde{a}_{\alpha+\beta},\widetilde{a}_{\alpha+\beta 2^s},\ldots,\widetilde{a}_{\alpha m^s+\beta (n-1)^s}]}}{({\alpha i^s+\beta j^s})-({\alpha +\beta})}
\end{equation*}
denotes $1^{st}$ and $(mn-1)^{th}$ Newton-divided difference of $\widetilde{a}_{\alpha+\beta}$ respectively, for Table \ref{table:1.0}, such that {$F[{\widetilde{a}_{\alpha i^s + \beta j^s}}] = \alpha i^s + \beta j^s$},  for $i=1,2,\ldots,m$ and $j=1,2,\ldots,n$.
\begin{table}[H]
\centering
\caption{An ordered arrangements of the elements of the given $({a_{ij}}) \in \mathcal{R}^{m \times n}$.}
{\begin{tabular}{|c|c|c|c|c|c|c|c|c|c|c|c|c|c|}
\hline
$\widetilde{a}_{\alpha+\beta}$ & ${\widetilde{a}_{\alpha + \beta 2^s}}$ & ... & ${\widetilde{a}_{\alpha + \beta n^s}}$ & ${\widetilde{a}_{\alpha 2^s + \beta}}$ & ... & ... & ${\widetilde{a}_{\alpha m^s + \beta }}$ & ${\widetilde{a}_{\alpha m^s + \beta 2^s}}$ & ... & ${\widetilde{a}_{\alpha m^s + \beta n^s}}$ 
 \\ 
 \hline
\end{tabular}}
\\where ${\widetilde{a}_{\alpha i^s + \beta j^s}}=a_{ij}$, for $i=1,2,\ldots,m$ and $j=1,2,\ldots,n.$
\label{table:1.0}
\end{table}
\end{cor}

Next, the particular choices $(\alpha,\beta) = (n,1)$ and $(\alpha,\beta) = (1,m)$ for $s=1$, lead to get the following results.

\begin{cor}\label{corollary:1.1}
For every $({a_{ij}}) \in \mathcal{R}^{m \times n}$, there exists a unique $p \in {_{n}^{1}\Pi_{mn-1}^{2}}$ that satisfies the DPPM \eqref{DPPM} is given as 
    \begin{equation}\label{pol:1}
	p(x,y)= \sum_{k=0}^{mn-1}{{{(nx+y)-(n+1)}\choose{k}}\Delta ^k {\widetilde{a}_{{1}}}},
	\end{equation}
where $\Delta^k {\widetilde{a}_{{1}}}$ is $k$th, $k=1,2,\ldots,mn-1$ Newton-forward difference of ${\widetilde{a}_{{1}}}$ for Table \ref{table:1}.
\begin{table}[H]
\centering
\caption{An ordered arrangements of the elements of the given $({a_{ij}}) \in \mathcal{R}^{m \times n}$.}
{\begin{tabular}{|c|c|c|c|c|c|c|c|c|c|c|c|c|c|}
\hline
 $\widetilde{a}_{1}$ & $\widetilde{a}_{2}$ & $\ldots$ & $\widetilde{a}_{n}$ & $\widetilde{a}_{n+1}$ & ... & ...  & $\widetilde{a}_{(m-1)n+1}$ & $\widetilde{a}_{(m-1)n+2}$ & ... & $\widetilde{a}_{mn}$ 
 \\ \hline
\end{tabular}}
\\where ${\widetilde{a}_{{(i-1)n+j}}}=a_{ij}$, for $i=1,2,\ldots,m$ and $j=1,2,\ldots,n$.
\label{table:1}
\end{table}
\end{cor}

\begin{proof}
Consider the injective transformation $\bar{\phi_1}:\mathcal{X} \rightarrow \mathcal{R}$ given by
\begin{equation} \label{pf:1}
		\bar{\phi_1}(i,j)=(i-1)n+j, \; \textrm{for} \; i=1,2,\ldots,m \textrm{ and } j=1,2,\ldots,n.
	\end{equation}
Using the transformation \eqref{pf:1}, the set of bivariate nodes $\mathcal{X}$ can be transformed into a sequence of length $mn$ and in an ordered arrangement can be written as 
\begin{flalign}\label{nodes:1.1} 
        \bar{\phi_1}{(\mathcal{X})} = \{\bar{\phi_1}(i,1),\bar{\phi_1}(i,2),...,\bar{\phi_1}(i,n): i = 1,2,...,m \}.
\end{flalign}
Here, the consecutive nodes in $\bar{\phi_1}{(\mathcal{X})}$ (in the given order) are equidistant with step size 1. Thus, using \eqref{pf:1}, the DPPM \eqref{DPPM} can be transformed into a univariate polynomial interpolation problem of the form
\begin{eqnarray}\label{DPPM:modified2}
            p(i,j) = \widetilde{p}(\bar{\phi_1}(i,j)) = \widetilde{p}((i-1)n + \beta j) = {a_{ij}} \textrm{ for all } (i,j) \in \mathcal{X}.
\end{eqnarray}
Therefore, on applying univariate Newton-forward difference interpolation formula for the interpolation problem \eqref{DPPM:modified2} with respect to $\mathcal{X}$, we get
\begin{equation}\label{ref:2}
	\begin{aligned}
	p(x,y) = {} & \widetilde{p}(\bar{\phi_1}(x,y)) = \widetilde{p}((x-1)n+y) \\ 
	= {} & {\widetilde{a}_{{1}}} + \frac{((nx+y)-(n+1))}{1!}{\Delta {\widetilde{a}_{{1}}}} \\ & + \frac{((nx+y)-(n+1))((nx+y)-(n+2))}{2!}{\Delta ^2 {\widetilde{a}_{{1}}}} \\ &\hspace{0.2cm} + 
	\ldots +\frac {(\bar{\phi_1}(x,y)-\bar{\phi_1}(1,1))...(\bar{\phi_1}(x,y)-\bar{\phi_1}(m,n-1))}{(mn-1)!}{\Delta^{mn-1} {\widetilde{a}_{{1}}}}.
	\end{aligned}
\end{equation}
On combining \eqref{ref:2} and \eqref{pf:1}, we get \eqref{pol:1}. Here, $p \in {_{n}^{1}\Pi_{mn-1}^{2}}$ and $(\alpha,\beta)=(n,1)$ satisfy the condition \eqref{cond:nodes} for $s=1$. This completes the result.
\end{proof}

\begin{cor} \label{corollary:1.2}
For every $({a_{ij}}) \in \mathcal{R}^{m \times n}$, there exists a unique $p \in {_{1}^{m}\Pi_{mn-1}^{2}}$ that satisfies the DPPM \eqref{DPPM} is given as
\begin{equation}\label{pol:2}
	p(x,y)= \sum_{k=0}^{mn-1}{{{(x+my)-(m+1)}\choose{k}}\Delta ^k \widetilde{a}_{1}},
\end{equation}
where $\Delta^k {\widetilde{a}_{{1}}}$ is $k$th, $k=1,2,\ldots,mn-1$ Newton-forward difference of ${\widetilde{a}_{{1}}}$ for Table \ref{table:2}.

\begin{table}[H]
\centering
\caption{An ordered arrangements of the elements of the given $({a_{ij}}) \in \mathcal{R}^{m \times n}$.}
{\begin{tabular}{|c|c|c|c|c|c|c|c|c|c|c|c|c|c|}
\hline
$\widetilde{a}_{1}$ & $\widetilde{a}_{2}$ & $\ldots$ & $\widetilde{a}_{m}$ & $\widetilde{a}_{m+1}$ & ... & ...  & $\widetilde{a}_{m(n-1)+1}$ & $\widetilde{a}_{m(n-1)+2}$ & ... & $\widetilde{a}_{mn}$ \\ 
\hline
\end{tabular}}
\\where ${\widetilde{a}_{{i+(j-1)m}}}=a_{ij}$, for $i=1,2,\ldots,m$ and $j=1,2,\ldots,n$.
\label{table:2}
\end{table}
\end{cor}

\begin{proof}
The proof is similar to Corollary \ref{corollary:1.1}. Particularly, the use of the injective transformation $\hat{\phi_1}:\mathcal{X} \rightarrow \mathcal{R}$, given as
\begin{equation*} \label{pf:2}
		\hat{\phi_1}(i,j)=i+(j-1)m \; \textrm{for all} \; i=1,2,\ldots,m;\; j=1,2,\ldots,n,
\end{equation*}
in place of \eqref{pf:1} completes the result.
\end{proof}

\begin{rem}\label{extension remark}
In general, there are uncountably many injective transformations from $\mathcal{X}$ to $\mathcal{R}$ that can project all the bivariate nodes on the grid $\mathcal{X}$ onto a finite set of $mn$ univariate points such that all projections are distinct. Therefore, the use of distinct injective transformation from $\mathcal{X}$ to $\mathcal{R}$, in place of \eqref{pf:11}, can lead us to get several classes of $mn$-dimensional subspaces of $\Pi^2$, of various total degrees, in which the DPPM \eqref{DPPM} can always be correct. However, for every distinct choice of such injective transformation, the resulting polynomials or correct spaces for the DPPM \eqref{DPPM} need not be distinct. For instance, on taking $\bar{\bar{\phi_1}},\hat{\hat{\phi_1}}:\mathcal{X} \rightarrow \mathcal{R}$, $\bar{\bar{\phi_1}}(i,j)= (mn+n+1)-ni-j$ and $\hat{\hat{\phi_1}}(i,j)= (mn+m+1)-i-mj$, in place of \eqref{pf:11}, the obtained polynomials coincide with the polynomials \eqref{pol:1} and \eqref{pol:2} in ${_{n}^{1}\Pi_{mn-1}^{2}}$ and ${_{1}^{m}\Pi_{mn-1}^{2}}$, respectively. 
\end{rem}

\section{An Application in Lagrange Bivariate Polynomial Interpolation}\label{sec4}
To be sure, together with $\mathcal{P}_m^n$, the LBVPIP \eqref{lag:int} with respect to $\mathcal{X}$ is also correct in $_{\alpha}^{\beta}\Pi_{s(mn-1)}^{2}$ for uncountably many choices of the pair $(\alpha,\beta)$, for each $s \in \mathbb{N}$. Moreover, for a bijective transformation $\psi:\Theta \rightarrow \mathcal{X}$, such that 
\begin{equation}\label{cond:5}
    \psi(x_i,y_j)=(i,j) \;\textrm{ for all } (i,j) \in \mathcal{X},
\end{equation}
the LBVPIP \eqref{lag:int} with respect to $\Theta$ can always be transformed into the LBVPIP \eqref{lag:int} with respect to $\mathcal{X}$, or equivalently into the DPPM \eqref{DPPM}. Particularly, if
\begin{equation*}
    x_i=x_1+(i-1)h, \; i=1,2,\ldots,m \quad \textrm{and} \quad y_j=y_1+(j-1)k, \; j=1,2,\ldots,n,  
\end{equation*}
where $h$ and $k$ be some real-scalars, then there exist a bijective map $\psi:\Theta \rightarrow \mathcal{X}$, given as 
\begin{equation*}
    \psi(x,y)=\left(1+\frac{x-x_1}{h},1+\frac{y-y_1}{k}\right) \; \textrm{for all} \; (x,y) \in \Theta,
\end{equation*}
which satisfies the condition \eqref{cond:5}. In this case, if a polynomial $p(x,y) \in \mathcal{P}$ satisfy the DPPM \eqref{DPPM}, then the polynomial 
\begin{equation*}
    p\left(x_1+(x-1)h, y_1+(y-1)k\right) \in \mathcal{P}
\end{equation*}
must satisfy the LBVPIP \eqref{lag:int} with respect to $\Theta$, i.e., if $\mathcal{P}$ is a correct space for the DPPM \eqref{DPPM}, then $\mathcal{P}$ must be a correct space for the LBVPIP \eqref{lag:int} with respect to $\Theta$. In extension, if the pair $(\alpha, \beta)$ satisfy the condition 
\begin{equation*}
    \alpha x_{1}^{s} + \beta y_{1}^{s} \neq \alpha x_{2}^{s} + \beta y_{2}^{s} \; \textrm{for all} \; (x_1,x_1) \neq (x_2,x_2),  
\end{equation*}
where $(x_1,y_1),(x_2,x_2) \in \Theta$. Then, the space $_{\alpha}^{\beta}\Pi_{s(mn-1)}^{2}$ is correct space for the LBVPIP \eqref{lag:int} with respect to $\Theta$ for uncountably many choices of the pair $(\alpha,\beta)$, for each $s \in \mathbb{N}$. 

Again let $I=\left[\min(x_1,x_2,\ldots,x_m),\max(x_1,x_2,\ldots,x_m)\right]$,\newline  $J=\left[\min(y_1,y_2,\ldots,y_n),\max(y_1,y_2,\ldots,y_n)\right]$ and $f$ be a sufficiently smooth real-valued bivariate function. Then, for some $\xi,\xi',\xi'' \in I$, $\eta,\eta',\eta'' \in J$, using the standard univariate error formulas recursively based on the mean value theorem or equivalently Taylor's series expansion of $f$ about $(x_1,y_1)$, the remainder terms $R(x,y)$ and $R_s(x,y)$ for the interpolating polynomials in $\mathcal{P}_{m}^{n}$ and $_{\alpha}^{\beta}\Pi_{s(mn-1)}^{2}$ respectively, for the LBVPIP \eqref{lag:int} for $\Theta$ can be written as
\begin{equation*}
        \begin{aligned}
    R(x,y) = \frac{1}{m!}{\frac{\partial^{m}f(\xi,y)}{\partial x^{m}}}{{\prod_{i=1}^{m}}{\left(x-x_i\right)}} + & \frac{1}{n!}{\frac{\partial^{n}f(x,\eta)}{\partial y^{n}}}{{\prod_{j=1}^{n}}{\left(y-y_j\right)}}\\& -\frac{1}{m!n!}{\frac{\partial^{m+n}f(\xi',\eta')}{\partial x^{m} \partial y^{n}}}{{\prod_{i=1}^{m}}{\left(x-x_i\right)}}{{\prod_{j=1}^{n}}{\left(y-y_j\right)}}
        \end{aligned}
\end{equation*}
and
\begin{equation*}
    R_s(x,y) = \frac{1}{(mn)!}{\sum_{k=0}^{mn}{{mn \choose{k}} \frac{\partial^{mn}{f(\xi'',\eta'')}}{{\partial x}^{mn-k} {{\partial y}^{k}}}}}{{\prod_{(i,j) \in \mathcal{X}}}{\left(\alpha (x^s-x_i^s) + \beta (y^s-y_j^s)\right)}}.
\end{equation*}
However, some limitations in the interpolation spaces $_{\alpha}^{\beta}\Pi_{s(mn-1)}^{2}$ for the LBVPIP \eqref{lag:int} with respect to $\Theta$ can be listed as follows:
\begin{itemize}
    \item The degree of the interpolating polynomials is rather high from the computational point of view.
    \item The interpolants are all constant in the direction $(-\beta,\alpha)$, an unfortunate bias.
\end{itemize} 

\section{Numerical Experiments}\label{sec5}
Let $P_{A} \in \mathcal{P}_{m}^{n}$, $p_{A} \in {_{n}^{1}\Pi_{mn-1}^{2}}$, and ${q}_{A} \in {_{1}^{m}\Pi_{mn-1}^{2}}$ respectively, represent the unique polynomials that satisfy the DPPM \eqref{DPPM} for $A \in \mathcal{R}^{m \times n}$. We use these notations with described meanings frequently throughout this section. 

\textbf{Example 1.} 
Using the result of Corollary \ref{cor:construction by lagrange} or Corollary \ref{cor:construction by newton forward}, the map $\mathcal{D}_P:\mathcal{R}^{2\times2} \rightarrow \mathcal{P}_{2}^{2}$, as defined as \eqref{dp isomorphism}, is given as
\begin{equation}\label{isomorphism:3}
    \begin{aligned}
	    \mathcal{D}_P{\left[\begin{pmatrix}		
	    a & b\\
		c & d
		\end{pmatrix} \right]} = (4a-2b-2c+&d) - (2a-b-2c+d)x \\ 
		& - (2a-2b-c+d)y + (a-b-c+d)xy. &&
		\end{aligned}
\end{equation}
For $A,B \in  \mathcal{R}^{2 \times 2}$ and non-zero scalar $\lambda$, the map \eqref{isomorphism:3} implies that
\begin{equation*}
    {\mathcal{D}_P}{(\lambda A+B)} = \lambda{ \mathcal{D}_P}(A)+{ \mathcal{D}_P}(B),
\end{equation*}
i.e., the map \eqref{isomorphism:3} is linear. Again, since  
$\mathcal{B}_1=\left\{ \begin{pmatrix}		
	    1 & 0\\
		0 & 0
		\end{pmatrix}, \begin{pmatrix}		
	    0 & 1\\
		0 & 0
		\end{pmatrix}, \begin{pmatrix}		
	    0 & 0\\
		1 & 0
		\end{pmatrix}, \begin{pmatrix}		
	    0 & 0\\
		0 & 1
		\end{pmatrix}\right\}$ and $\mathcal{B}_2={\{1,x,y,xy\}}$ be the sets of standard bases of $ \mathcal{R}^{2 \times 2}$ and $\mathcal{P}_{2}^{2}$, respectively. Then, using \eqref{isomorphism:3}, we get
\begin{equation*}
    \mathcal{D}_P{\left[\begin{pmatrix}		
	    1 & 0\\
		0 & 0
		\end{pmatrix} \right]} = 4-2x-2y+xy,
\end{equation*}		
\begin{equation*}
    \mathcal{D}_P{\left[\begin{pmatrix}		
	    0 & 1\\
		0 & 0
		\end{pmatrix} \right]} = -2+x+2y-xy,
\end{equation*}
\begin{equation*}
    \mathcal{D}_P{\left[\begin{pmatrix}		
	    0 & 0\\
		1 & 0
		\end{pmatrix} \right]} =  -2+2x+y-xy,
\end{equation*}
and
\begin{equation*}
    \mathcal{D}_P{\left[\begin{pmatrix}		
	    0 & 0\\
		0 & 1
		\end{pmatrix} \right]} = 1-x-y+xy.
\end{equation*}
Therefore, the co-ordinate matrix of the map \eqref{isomorphism:3}, with respect to bases $\mathcal{B}_1$ and $\mathcal{B}_2$, is written as
\begin{equation}\label{co-ordinate matrix 1}
	   [ \mathcal{D}_P]_{\mathcal{B}_1}^{\mathcal{B}_2} = \begin{pmatrix}		
	    4 & -2 & -2 & 1\\
		-2 & 1 & 2 & -1\\
		-2 & 2 & 1 & -1\\
		1 & -1 & -1 & 1
		\end{pmatrix}.
\end{equation}
This implies, $det\left([ \mathcal{D}_P]_{\mathcal{B}_1}^{\mathcal{B}_2} \right)=-1 \neq 0$, i.e., $ker( \mathcal{D}_P)=\left\{0\right\}$, i.e., the map \eqref{isomorphism:3} is invertible. This verifies that the map \eqref{isomorphism:3} is an isomorphism that satisfies the DPPM \eqref{DPPM}. Moreover, using \eqref{inverse isomorphic map2.1}, the map $\mathcal{D}_P^{-1}:\mathcal{P}_{2}^{2} \rightarrow  \mathcal{R}^{2 \times 2}$ is given as \\
$\mathcal{D}_P^{-1}{\left(\lambda_{0,0}+\lambda_{1,0}x+\lambda_{0,1}y+\lambda_{1,1}xy\right)}$
\begin{equation}\label{inverse map 3}
	   = \begin{pmatrix}		
	    \lambda_{0,0}+\lambda_{1,0}+\lambda_{0,1}+\lambda_{1,1} & \lambda_{0,0}+\lambda_{1,0}+2\lambda_{0,1}+2\lambda_{1,1}\\
		\lambda_{0,0}+2\lambda_{1,0}+\lambda_{0,1}+2\lambda_{1,1} & \lambda_{0,0}+2\lambda_{1,0}+2\lambda_{0,1}+4\lambda_{1,1}
		\end{pmatrix}, 
\end{equation}
for all $\lambda_{0,0}+\lambda_{1,0}x+\lambda_{0,1}y+\lambda_{1,1}xy \in \mathcal{P}_{2}^{2}$. Then, using \eqref{inverse map 3}, we get
\begin{equation*}
    \mathcal{D}_P^{-1}(1) = 1\begin{pmatrix}		
	    1 & 0\\
		0 & 0
		\end{pmatrix} + 1\begin{pmatrix}		
	    0 & 1\\
		0 & 0
		\end{pmatrix} + 1\begin{pmatrix}		
	    0 & 0\\
		1 & 0
		\end{pmatrix} + 1\begin{pmatrix}		
	    0 & 0\\
		0 & 1
		\end{pmatrix},
\end{equation*}
\begin{equation*}
    \mathcal{D}_P^{-1}(x) = 1\begin{pmatrix}		
	    1 & 0\\
		0 & 0
		\end{pmatrix} + 1\begin{pmatrix}		
	    0 & 1\\
		0 & 0
		\end{pmatrix} + 2\begin{pmatrix}		
	    0 & 0\\
		1 & 0
		\end{pmatrix} + 2\begin{pmatrix}		
	    0 & 0\\
		0 & 1
		\end{pmatrix},
\end{equation*}
\begin{equation*}
    \mathcal{D}_P^{-1}(y) = 1\begin{pmatrix}		
	    1 & 0\\
		0 & 0
		\end{pmatrix} + 2\begin{pmatrix}		
	    0 & 1\\
		0 & 0
		\end{pmatrix} + 1\begin{pmatrix}		
	    0 & 0\\
		1 & 0
		\end{pmatrix} + 2\begin{pmatrix}		
	    0 & 0\\
		0 & 1
		\end{pmatrix},
\end{equation*}
and
\begin{equation*}
    \mathcal{D}_P^{-1}(xy) = 1\begin{pmatrix}		
	    1 & 0\\
		0 & 0
		\end{pmatrix} + 2\begin{pmatrix}		
	    0 & 1\\
		0 & 0
		\end{pmatrix} + 2\begin{pmatrix}		
	    0 & 0\\
		1 & 0
		\end{pmatrix} + 4\begin{pmatrix}		
	    0 & 0\\
		0 & 1
		\end{pmatrix}.
\end{equation*}
Therefore, the co-ordinate matrix of the map \eqref{inverse map 3}, with respect to bases $\mathcal{B}_2$ and $\mathcal{B}_1$, is written as
\begin{equation}\label{co-ordinate matrix 2}
	   [\mathcal{D}_P^{-1}]_{\mathcal{B}_2}^{\mathcal{B}_1} = \begin{pmatrix}		
	    1 & 1 & 1 & 1\\
		1 & 1 & 2 & 2\\
		1 & 2 & 1 & 2\\
		1 & 2 & 2 & 4
		\end{pmatrix},
\end{equation}
such that $det\left([\mathcal{D}_P^{-1}]_{\mathcal{B}_2}^{\mathcal{B}_1} \right)=-1 \neq 0$. Again, using \eqref{co-ordinate matrix 1} and \eqref{co-ordinate matrix 2}, we get
\begin{equation*}
    [\mathcal{D}_P]_{\mathcal{B}_1}^{\mathcal{B}_2}.[\mathcal{D}_P^{-1}]_{\mathcal{B}_2}^{\mathcal{B}_1}=\begin{pmatrix}		
	    1 & 0 & 0 & 0\\
	    0 & 1 & 0 & 0\\
	    0 & 0 & 1 & 0\\
	    0 & 0 & 0 & 1
		\end{pmatrix}=[\mathcal{D}_P^{-1}]_{\mathcal{B}_2}^{\mathcal{B}_1}.[ \mathcal{D}_P]_{\mathcal{B}_1}^{\mathcal{B}_2}.
\end{equation*}
This verifies that the map \eqref{inverse map 3} is the inverse map of \eqref{isomorphism:3}. 
In a similar manner, using the results of Corollary \ref{corollary:1.1} and  Corollary \ref{corollary:1.2}, the isomorphisms $\mathcal{D}_p:\mathcal{R}^{2\times2} \rightarrow \;_{2}^{1}\Pi_{3}^{2}$ and $\mathcal{D}_q:\mathcal{R}^{2\times2} \rightarrow \;_{1}^{2}\Pi_{3}^{2}$ respectively, as defined as \eqref{dp isomorphism}, are given as
\begin{equation}\label{isomorphism:1}
    \begin{aligned}
	   \mathcal{D}_p{\left[\begin{pmatrix}		
	    a & b\\
		c & d
	\end{pmatrix} \right]} = {} & (20a-45b+36c-10d) + \frac{(-74a+189b-162c+47d)}{6}(2x+y) \\
	& + \frac{(5a-14b+13c-4d)}{2}(2x+y)^2 + \frac{(-a+3b-3c+d)}{6}(2x+y)^3 
		\end{aligned}
\end{equation}
and
\begin{equation}\label{isomorphism:2}
    \begin{aligned}
	   \mathcal{D}_q{\left[\begin{pmatrix}		
	    a & b\\
		c & d
	\end{pmatrix} \right]} = {} & (20a+36b-45c-10d) + \frac{(-74a-162b+189c+47d)}{6}(x+2y) \\
	& + \frac{(5a+13b-14c-4d)}{2}(x+2y)^2 + \frac{(-a-3b+3c+d)}{6}(x+2y)^3. 
		\end{aligned}
\end{equation}
Also, using Remark \ref{inverse: isomorphic map def}, the inverse linear maps $\mathcal{D}_{p}^{-1}:\;_{2}^{1}\Pi_{3}^{2}\rightarrow \mathcal{R}^{2 \times 2}$ and $\mathcal{D}_{q}^{-1}:\;_{1}^{2}\Pi_{3}^{2}\rightarrow \mathcal{R}^{2 \times 2}$ of the maps \eqref{isomorphism:1} and \eqref{isomorphism:2} respectively, are given as 
\begin{equation}\label{inverse DP map:2.1}
	   \mathcal{D}_p^{-1}{\left(p+q u+r u^2+s u^3 \right)} = \begin{pmatrix}		
	    p+3q+9r+27s & p+4q+16r+64s\\
		p+5q+25r+125s & p+6q+36r+216s
		\end{pmatrix} \nonumber
\end{equation}
and
\begin{equation}\label{inverse DP map:2.2}
	   \mathcal{D}_q^{-1}{\left(p+q v+r v^2+s v^3 \right)} = \begin{pmatrix}		
	    p+3q+9r+27s & p+5q+25r+125s\\
		p+4q+16r+64s & p+6q+36r+216s
		\end{pmatrix},   \nonumber
\end{equation}
where $u=2x+y$, $v=x+2y$ and $p$, $q$, $r$, $s$ are some scalars. Clearly, $p+q u+r u^2+s u^3 \in {}_{2}^{1}\Pi_{3}^{2}$ and $p+q v+r v^2+s v^3 \in {}_{1}^{2}\Pi_{3}^{2}$.\\

\textbf{Example 2.} 
Let $\delta \in \mathcal{R}^{1 \times 3}$ and $\vartheta \in \mathcal{R}^{2 \times 2}$ be two given matrices, such as 
\begin{equation}
    \delta = \begin{pmatrix}			
		1 & -1 & -2
		    \end{pmatrix} \textrm{ and } \vartheta = \begin{pmatrix}     
		-15 & 36\\
		-1 & 96
		    \end{pmatrix}.
\end{equation}
The DPPM \eqref{DPPM} for $\delta \in \mathcal{R}^{1 \times 3}$ has no solution in 3-dimensional subspace $\Pi_1^2 \subset \Pi^2$ and the DPPM \eqref{DPPM} for $\vartheta \in \mathcal{R}^{2 \times 2}$ does not hold the necessary condition of given data points in $(k+1)(k+2)/2$-dimensional subspace $\Pi_k^2 \subset \Pi^2$ for any $k \in \mathbb{N}$. In other words, the DPPM \eqref{DPPM} for $\delta \in \mathcal{R}^{1 \times 3}$ or $\vartheta \in \mathcal{R}^{2 \times 2}$ does not possess any solution in $\Pi_k^2$ for any $k \in \mathbb{N}$. Now, using the results of Corollary \ref{cor:construction by lagrange} or \ref{cor:construction by newton forward}, Corollary \ref{corollary:1.1}, and Corollary \ref{corollary:1.2} respectively, we get 
\begin{equation}
    {P}_{\delta}(x,y) = {~} \frac{1}{2}y^2-\frac{7}{2}y+4 \in \mathcal{P}_{1}^{3}, \label{P:delta}
\end{equation}	
\begin{equation}
    {P}_\vartheta(x,y) = {~} 46xy-32x+5y-34 \in \mathcal{P}_{2}^{2}, \label{P:vartheta}
\end{equation}	
\begin{equation}
    p_{\delta}(x,y) = {~} \frac{1}{2}(3x+y)^2-\frac{13}{2}(3x+y)+19 \in {_{3}^{1}\Pi_{2}^{2}}, \label{p:delta}
\end{equation}	
\begin{equation}
    p_\vartheta(x,y) = {~} 37(2x+y)^3-488(2x+y)^2+2098(2x+y)-2916 \in {_{2}^{1}\Pi_{3}^{2}},  \label{p:vartheta}
\end{equation}	
\begin{equation}
    {q}_{\delta}(x,y) = {~} \frac{1}{2}(x+y)^2-\frac{9}{2}(x+y)+8 \in {_{1}^{1}\Pi_{2}^{2}}, \label{q:delta}
\end{equation}	
and
\begin{equation}
    {q}_\vartheta(x,y) = {~} \frac{23}{2}(x+2y)^2-\frac{133}{2}(x+2y)+81 \in {_{1}^{2}\Pi_{3}^{2}}.  \label{q:vartheta}
\end{equation}	
\begin{figure}[H]
    \centering
    \begin{subfigure}[b]{0.4\textwidth}
        \includegraphics[width=\textwidth]{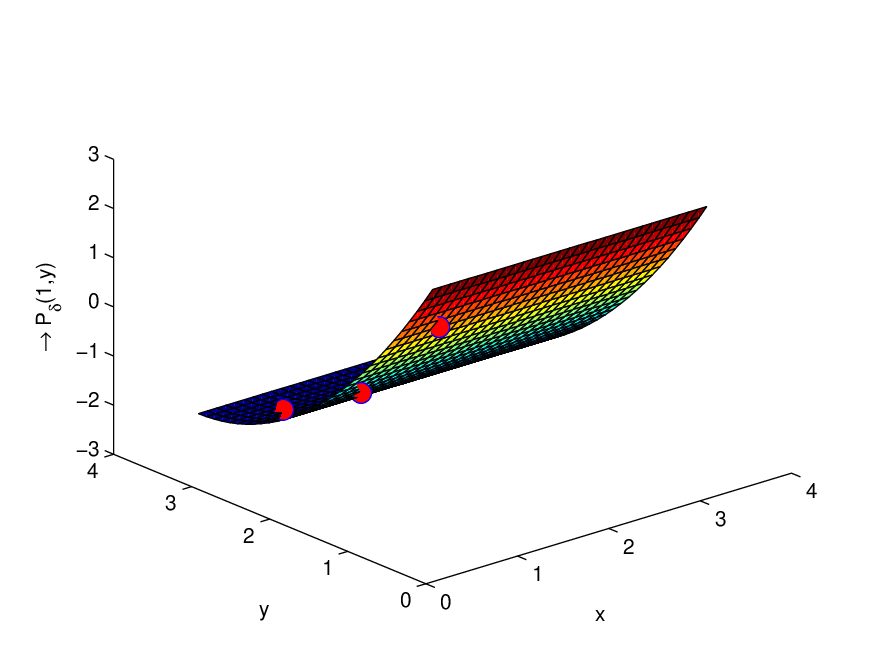}
        \caption{$P_{\delta} \in \mathcal{P}_{2}^{2}$}
        \label{fig:delta3}
    \end{subfigure}
~     \begin{subfigure}[b]{0.4\textwidth}
        \includegraphics[width=\textwidth]{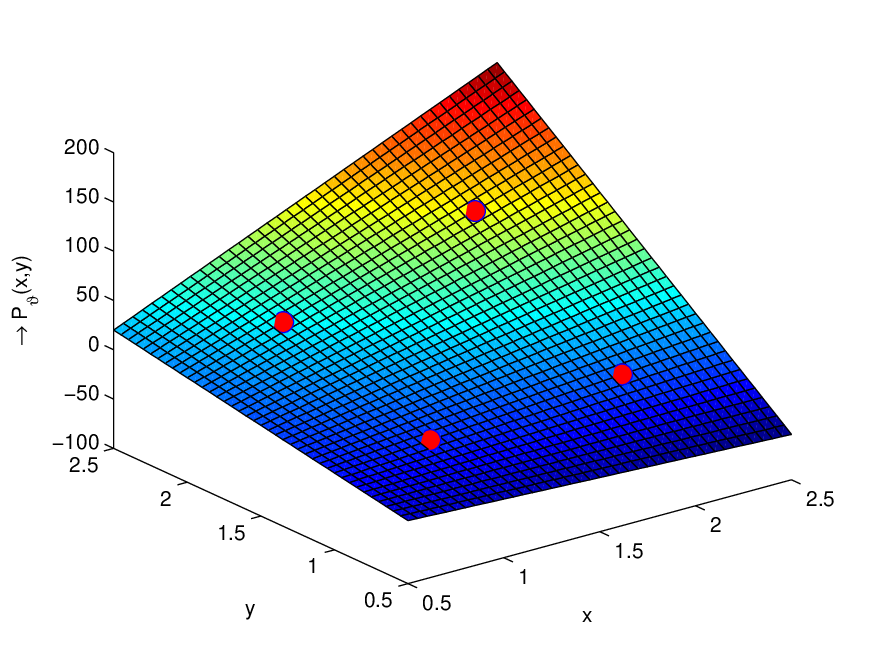}
        \caption{${P}_{\vartheta} \in \mathcal{P}_{2}^{2}$}
        \label{fig:vartheta3}
    \end{subfigure}
\hfill
    \begin{subfigure}[b]{0.4\textwidth}
        \includegraphics[width=\textwidth]{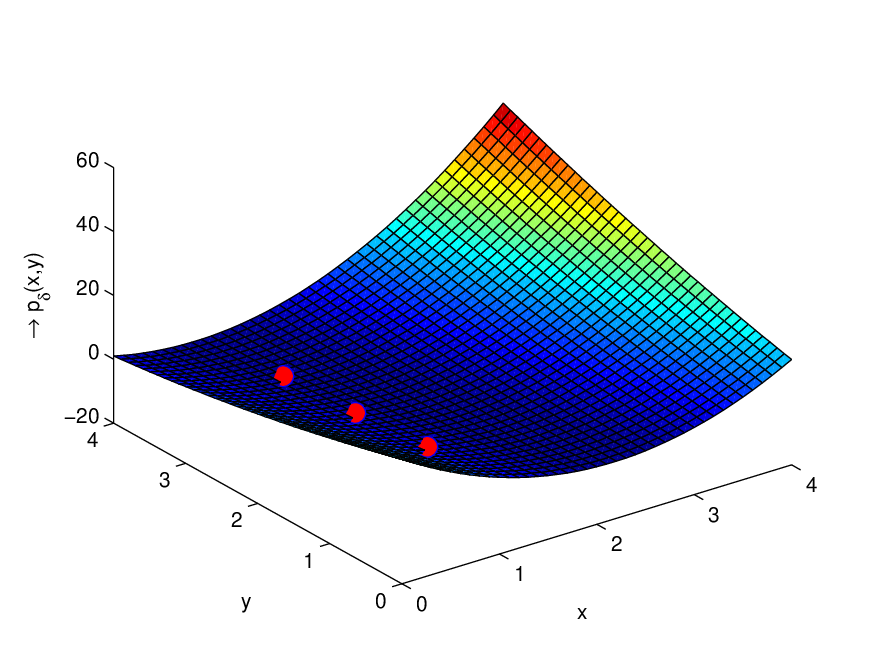}
        \caption{$p_{\delta} \in {_{3}^{1}\Pi_{1,3}^{2}}$}
        \label{fig:delta1}
    \end{subfigure}
~   \begin{subfigure}[b]{0.4\textwidth}
        \includegraphics[width=\textwidth]{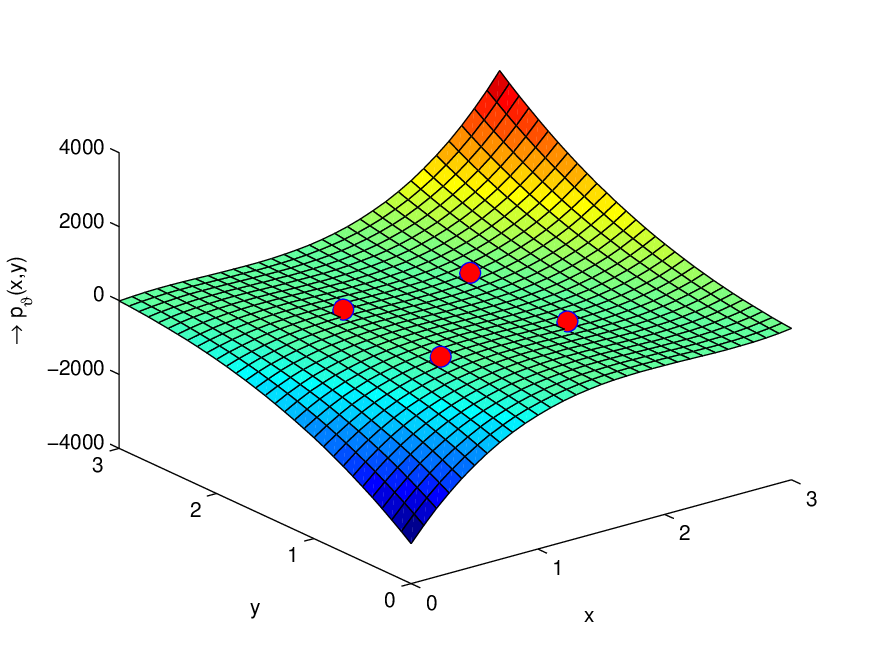}
        \caption{${p}_{\vartheta} \in {_{1}^{1}\Pi_{1,3}^{2}}$}
        \label{fig:vartheta1}
    \end{subfigure}
\hfill
    \begin{subfigure}[b]{0.4\textwidth}
        \includegraphics[width=\textwidth]{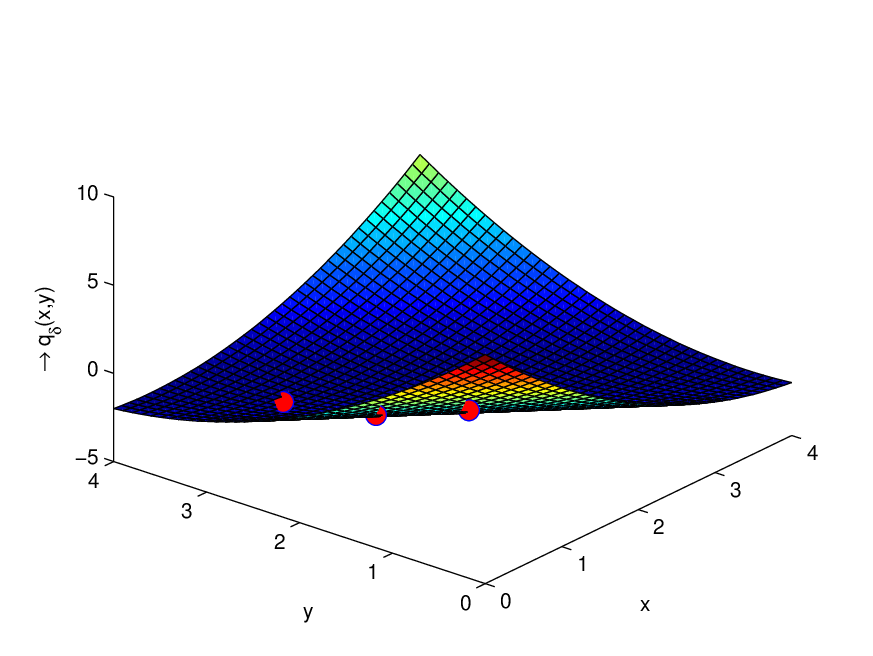}
        \caption{$q_{\delta} \in {_{2}^{1}\Pi_{3}^{2}}$}
        \label{fig:delta2}
    \end{subfigure}
~   \begin{subfigure}[b]{0.4\textwidth}
        \includegraphics[width=\textwidth]{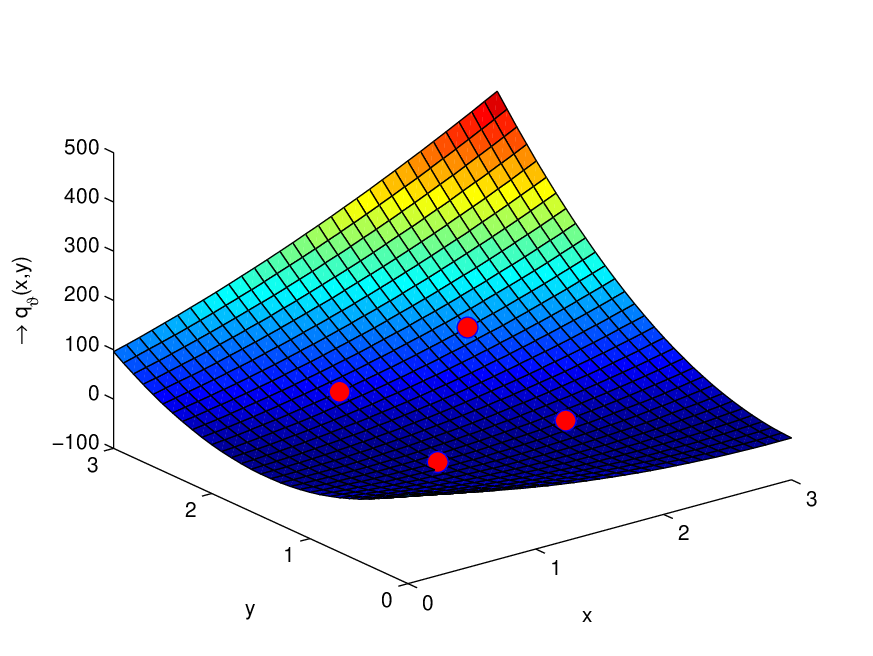}
        \caption{${q}_{\vartheta} \in {_{1}^{2}\Pi_{3}^{2}}$}
        \label{fig:vartheta2}
    \end{subfigure} 
    \caption{Surface diagrams of the polynomials $P_{\delta} \in \mathcal{P}_{1}^{3}$, ${P}_{\vartheta} \in \mathcal{P}_{2}^{2}$, $p_{\delta} \in {_{3}^{1}\Pi_{2}^{2}}$, $p_{\vartheta} \in {_{2}^{1}\Pi_{3}^{2}}$, ${q}_{\delta} \in {_{1}^{1}\Pi_{2}^{2}}$ and ${q}_{\vartheta} \in {_{1}^{2}\Pi_{3}^{2}}$ indicating all the data points of the given matrices $\delta \in \mathcal{R}^{1 \times 3}$ and $\vartheta \in \mathcal{R}^{2 \times 2}$ respectively.}\label{surface:diagram2}
\end{figure}

\textbf{Example 3.} 
Let us consider a LBVPIP \eqref{lag:int} with respect to the set of rectangular bivariate nodes $\bar{\mathcal{X}}=\{(1,1),(1,2)(2,1),(2,2)\}$, concerning the data points for the real-valued function 
\begin{equation}\label{function:example}
    f(x,y)= 3{x^2}+2{y^3}-{x^2}y+2x{y^2}+x-y+10.
\end{equation}
Without loss of generality, the LBVPIP \eqref{lag:int} with respect to $\bar{\mathcal{X}}$ can be converted into the DPPM \eqref{DPPM} for $\kappa = \left(f(i,j)\right) \in \mathcal{R}^{2 \times 2}$. Therefore, using the results of Corollary \ref{cor:construction by lagrange} or \ref{cor:construction by newton forward} and Corollary \ref{corollary:1.2}, we get 
\begin{equation}\label{pol:ex2.1}
    P_{\kappa}(x,y)  = {~} 3xy + 6x + 15y - 8
\end{equation}
and
\begin{equation}\label{pol:ex2.2}
    {q}_{\kappa}(x,y)  = {~} \frac{1}{2}{(x+2y)^3} - 6(x+2y)^2 + \frac{65}{2}{(x+2y)} -41.
\end{equation}
Now, using \eqref{function:example} and \eqref{pol:ex2.1}, the absolute error of the interpolating polynomial $P_{\kappa} \in \mathcal{P}_{2}^{2}$ can be given as
\begin{equation}\label{error:2}
    E_{P_{\kappa} \in {~}\mathcal{P}_{2}^{2}} = |{f(x,y)-{P}_{\kappa}(x,y)}|, \quad (x,y) \in [1,2]\times[1,2].
\end{equation}
Also, using \eqref{function:example} and \eqref{pol:ex2.2}, the absolute error of the interpolating polynomial ${q}_{\kappa} \in {_{1}^{2}\Pi_{3}^{2}}$ can be given as 
\begin{equation}\label{error:1}
    E_{{q}_{\kappa} \in {~}{_{1}^{2}\Pi_{3}^{2}}} = |{f(x,y)-{q}_{\kappa}(x,y)}|, \quad (x,y) \in [1,2]\times[1,2].
\end{equation}
Next, in Table \ref{table:3}, we calculate the difference of absolute errors of the interpolating polynomials in ${q}_{\kappa} \in {_{1}^{2}\Pi_{3}^{2}}$ and $P_{\kappa} \in \mathcal{P}_{2}^{2}$ respectively, for the LBVPIP \eqref{lag:int} for $\bar{\mathcal{X}}$ at the equally spaced mesh-grid points $(x,y) \in [1,2]\times[1,2]$ with step size 0.1, up to 4 decimal places. 
\begin{table}[H]
\centering
\caption{The difference of absolute errors for ${q}_{\kappa} \in {_{1}^{2}\Pi_{3}^{2}}$ and $P_{\kappa} \in \mathcal{P}_{2}^{2}$ in $[1,2]\times[1,2]$.} 
{\begin{tabular}{|c|c|c|c|c|c|c|}
\hline
\multicolumn{7}{|c|}{The difference of absolute errors}
\\
\multicolumn{7}{|c|}{$E_{{q_{\kappa}} \in {_{1}^{2}\Pi_{3}^{2}}} -  E_{P_{\kappa} \in \mathcal{P}_{2}^{2}}$ =  $|{f(x,y)-{q}_{\kappa}(x,y)}| - |{f(x,y)-P_{\kappa}(x,y)}|$}
\\
\multicolumn{7}{|c|}{at the equally spaced mesh-grid points $(x,y) \in [1,2]\times[1,2]$}\\
\multicolumn{7}{|c|}{with step size 0.1 by 0.2, up to 4 decimal places.}
\\
\hline
 $\downarrow \textbf{y}$ $\rightarrow \textbf{x}$ & \bf{1.0} &  \bf{1.2}  & \bf{1.4} &  \bf{1.6}  & \bf{1.8} &  \bf{2.0}
 \\ 
\hline
 \bf{1.0} & 0 &  0.1439  & 0.1920 &  0.1679  & 0.0959 &  0
 \\ 
\hline
\bf{1.1} & 0.1439 &  0.1320  & 0.0479  &   -0.0839 & -0.2400  & -0.3959
 \\ 
\hline
\bf{1.2} & 0.1920  & 0.0479 & -0.1440 &  -0.3599 &-0.5759 & -0.7680
\\ 
\hline
\bf{1.3} & 0.1679 & -0.0840 & -0.3599 &   -0.6359 & -0.8880 & -1.0919
\\ 
\hline
\bf{1.4} &   0.0959 & -0.2399 & -0.5760 &   -0.8880 & -1.1519 & -1.3440
\\ 
\hline
\bf{1.5} & 0 & -0.3959  & -0.7680 &   -1.0919 & -1.3440 &  -1.5000
\\ 
\hline
\bf{1.6} &  -0.0959 & -0.5280 & -0.9119 &   -1.2239 & -1.4399 & -1.5360
\\ 
\hline
\bf{1.7} &  -0.1680 & -0.6119 & -0.9840 &  -1.2600 & -1.4160 & -1.4280
\\ 
\hline
\bf{1.8} &  -0.1919 & -0.6240  & -0.9600 &    -1.1760 &  -1.2480 &  -1.1519
\\ 
\hline
\bf{1.9} &  -0.1440  & -0.5399  & -0.8160  &   -0.9480  & -0.9119  & -0.6840
\\ 
\hline
\bf{2.0} & 0  & 0.0159   &  0.0480    &   0.0719   &  0.0640    & 0
\\ 
\hline
\end{tabular}}
\label{table:3}
\end{table}
From Table \ref{table:3}, we can observe that the absolute error of the interpolating polynomial ${q_{\kappa}}\in {_{1}^{2}\Pi_{3}^{2}}$ is moreover less (the values indicating `-' ve sign) in comparison to the absolute error the interpolating polynomial $P_{\kappa}\in \mathcal{P}_{2}^{2}$, at most of the interpolating points in the grid $[1,2]\times[1,2]$. 

\section{Concluding Remarks}\label{sec6}
For a given matrix $({a_{ij}}) \in \mathcal{R}^{m \times n}$ and an $mn$-dimensional subspace $\mathcal{P} \subset \Pi^2$, the article proposes a bivariate polynomial problem for matrices, namely DPPM \eqref{DPPM}. As an advanced application, the DPPM \eqref{DPPM} provides a sufficient condition for the map \eqref{dp isomorphism} to be an isomorphism. In the process of solving, it is established that the DPPM \eqref{DPPM} and the LBVPIP \eqref{lag:int} for $\mathcal{X}$ are equivalent. Theorem \ref{sufficient condition} proves that there always exists an isomorphism, as defined as \eqref{dp isomorphism}, if $\mathcal{P}$ is a correct space for the DPPM \eqref{DPPM}. As well, Remark \ref{inverse: isomorphic map def} provides the associated inverse linear map of the map \eqref{dp isomorphism} and Corollary \ref{special property} proposes an important consequence to the DPPM \eqref{DPPM}. Further, Theorem \ref{thm:dp} proves the existence and uniqueness of the solution of the DPPM \eqref{DPPM} in $\mathcal{P}_{m}^{n}$. For this, the standard tensor product approach of univariate polynomials is used for the DPPM \eqref{DPPM}. Also, Corollary \ref{cor:construction by lagrange} and Corollary \ref{cor:construction by newton forward} provide two formulas to construct the required unique polynomial in $\mathcal{P}_{m}^{n}$ that satisfies the DPPM\eqref{DPPM}, which involve univariate Newton-Lagrange's and Newton-forward differences of data values, respectively. 

Thereafter, Definition \ref{new:subspace} and Remark \ref{basis2} introduce an $mn$-dimensional subspace of $\Pi^2_{s(mn-1)}$, involving two real parameters $\alpha$ and $\beta$, denoted as ${}_{\alpha}^{\beta}\Pi_{s(mn-1)}^{2}$. Theorem \ref{thm:gen1.1} proves that there always exists a polynomial in ${}_{\alpha}^{\beta}\Pi_{s(mn-1)}^{2}$ that uniquely satisfies the DPPM \eqref{DPPM} for uncountably many choices of the pair $(\alpha,\beta)$. This establishes the existence of a new class of $mn$-dimensional subspaces of $\Pi^2_{s(mn-1)}$, of various total degrees, for each $s \in \mathbb{N}$, in which the DPPM \eqref{DPPM} always possesses a unique solution. For this, a transformation $\phi_s:\mathcal{X} \rightarrow \mathbb{R}$,$(i,j) \rightarrow {\alpha i^s + \beta j^s}$ and $z = \alpha x^s + \beta y^s$, for some $\alpha > 0$ and $\beta > 0$, is used such that $\phi_s$ is injective. In other words, all points of $\mathcal{X}$ are projected to the curve orthogonal to $\alpha x^s +\beta y^s =0$, which are chosen so that all projections differ. As well, Corollary \ref{corollary:1.0} offers a common formula to construct the required unique polynomial that satisfies the DPPM \eqref{DPPM} in ${_{\alpha}^{\beta}\Pi_{s(mn-1)}^{2}}$, for each $s \in \mathbb{N}$. For computational simplicity, the particular choices of $(\alpha,\beta) = (n,1)$ and $(\alpha,\beta) = (1,m)$ with $s=1$ permit writing the interpolant in terms of the standard formulas involving forward differences of the data values for polynomial interpolation at consecutive integers, which are submitted in Corollary \ref{corollary:1.1} and Corollary \ref{corollary:1.2}, respectively. In extension, Remark \ref{extension remark} asserts that using distinct injective transformations from $\mathcal{X}$ to $\mathcal{R}$, several classes of $mn$-dimensional subspaces of $\Pi^2$, of various total degrees, can be obtained in which the DPPM \eqref{DPPM} can also be solved uniquely. 

In addition, section \ref{sec4} provides possible applicability of the results in the LBVPIP \eqref{lag:int} for $\mathcal{X}$ and the LBVPIP \eqref{lag:int} for $\Theta$, respectively. The remainder formulas are also submitted for the interpolating polynomials in $\mathcal{P}_{m}^{n}$ and $_{\alpha}^{\beta}\Pi_{s(mn-1)}^{2}$ respectively, concerning the LBVPIP \eqref{lag:int} for $\Theta$. Moreover, some limitations in the interpolation spaces $_{\alpha}^{\beta}\Pi_{s(mn-1)}^{2}$, $s \in \mathbb{N}$ have been discussed for the same. At last, Example 1, Example 2, and Example 3 validate, compare, and justify the theoretical findings and show the possible applicability of the results in linear algebra, numerical linear algebra, and interpolation theory. Particularly, Example 3 indicates that together with $\mathcal{P}_m^n$, the interpolation spaces $_{\alpha}^{\beta}\Pi_{s(mn-1)}^{2}$ can also be used to optimize the accuracy of the LBVPIP \eqref{lag:int} for the rectangular schemes of $m$ by $n$ bivariate nodes on $\mathcal{X}$, for some suitable choices of the pair $(\alpha,\beta)$ and $s \in \mathbb{N}$ that satisfy the condition \eqref{cond:nodes}. 

Overall, the article appears to be new results and a significant contribution to linear algebra and bivariate Lagrange polynomial interpolation, well-recognized and historically rich fields. In extension, some algebraic and geometric properties of bivariate polynomials can be constructed in $\mathcal{P}_{m}^{n}$ or $_{\alpha}^{\beta}\Pi_{s(mn-1)}^{2}$, corresponding to DPPM \eqref{DPPM} and the isomorphism \eqref{dp isomorphism}, which can preserve the algebraic and geometric structures of $\mathcal{R}^{m \times n}$. Such as algebra structure, ring structure, inner product structure, norm structure, metric-space structure, etc. The authors consider them as one of their future objectives.

\section{Material and Methods}\label{sec8}

The figures given in Example 2 of the section \ref{sec5} were obtained by using codes elaborated with software MATLAB \cite{21} without using any extra precision. All computations were performed on the HP Envy Leapmotion Touchsmart 17 inch 2014 edition Laptop (Intel(R) Core(TM) i7-4702MQ CPU @ 2.20GHz Processor with 8GB DDR3 SDRAM).\\

\textbf{Author contributions:}\\
\textit{Conceptualization $\&$ Writing-Original Draft:} Dharm Prakash Singh; \textit{Guidance $\&$ Modifications:} Amit Ujlayan; \textit{Discussions $\&$ Validations:} Bhim Sen Choudhary. \\

\textbf{Conflicts of Interest:} The authors declare no conflict of interest.


\bibliography{mybib}
\bibliographystyle{unsrt}

\end{document}